\documentclass{amsart}

\usepackage{amssymb,graphicx, color}
\usepackage{xypic}
\usepackage[all]{xy}
\usepackage{lscape}
\usepackage{longtable}

\newtheorem{thm}{Theorem}[section]

\newtheorem{conj}[thm]{Conjecture}
   
\theoremstyle{definition}
\newtheorem{defn}[thm]{Definition}

\newtheorem{rem}[thm]{Remark}

\newtheorem{defn-thm}[thm]{Definition--Theorem}  
\newtheorem{defn-lem}[thm]{Definition--Lemma}  

\theoremstyle{remark}

\renewcommand{\c}[0]{{\mathbb C}}

\newcommand{\C}{\mathbb{C}}

\def\loccoh#1.#2.#3.#4.{H^{#1}_{#2}(#3,#4)}

\DeclareMathAlphabet{\mathchanc}{OT1}{pzc}%
                                {m}{it}

\numberwithin{equation}{section}

\begin{document}
\bibliographystyle{amsalpha}

\title[A counterexample to a conjecture on simultaneous Waring identifiability]{A counterexample to a conjecture on simultaneous Waring identifiability}

\author[Elena Angelini]{Elena Angelini}
\email{elena.angelini@unisi.it}

\begin{abstract}
The new identifiable case appeared in \cite{AGMO}, together with the analysis on simultaneous identifiability of pairs of ternary forms recently developed in \cite{BG}, suggested the following conjecture towards a complete classification of all simultaneous Waring identifiable cases: for any $ d \geq 2 $, the general polynomial vectors consisting of $ d-1 $ ternary forms of degree $ d $ and a ternary form of degree $ d+1 $, with rank $ \frac{d^2+d+2}{2} $, are identifiable over $ \c $.\\
In this paper, by means of a computer-aided procedure inspired to the one described in \cite{AGMO}, we obtain that the case $ d = 4 $ contradicts the previous conjecture, admitting at least $ 36 $ complex simultaneous Waring decompositions (of length $ 11 $) instead of $ 1 $.   


\noindent \emph{Mathematics~Subject~Classification~(2020):}
14N07, 15A69, 14N05, 13P05, 65H10, 15A72, 14Q20.
\end{abstract}

\maketitle

\section{Introduction}\label{sec:intr}

The main purpose of this paper is to provide a contribution towards a complete classification of all \emph{simultaneous Waring identifiable} cases. In particular, we focus on general \emph{polynomial vectors}, which means vectors $ f= (f_{1}, \ldots, f_{r}) $ consisting of general homogeneous polynomials (for simplicity, forms) $ f_{j} \in \C[x_{0}, \ldots, x_{n}]_{d_{j}} $, with degrees $ d_{1} \leq \ldots \leq d_{r} $ and we are interested in the ways of expressing all the $ f_{i}$'s as a linear combination of powers (which explains the use of \emph{Waring}, in honor of the homonymous statement in Number Theory established in 1770) of the same (which justifies the use of the term \emph{simultaneous}) linear forms $ \ell_{i} \in \C[x_{0}, \ldots, x_{n}]_{1} $. Such expressions are thus called \emph{simultaneous Waring decompositions}. In this setting, two natural questions arise:
\begin{itemize}
\item[\bf{Q1}] what is the minimum number of linear forms necessary to express simultaneously $ f_{1}, \ldots, f_{r} $?
\item[\bf{Q2}] provided that we know the number introduced in Q1, which is called \emph{simultaneous Waring rank} of $ f= (f_{1}, \ldots, f_{r}) $ (or simply, rank) and is denoted by $ k $, does $ f $ admit a unique simultaneous Waring decomposition with $ k $ summands? 
\end{itemize}
If Q2 has positive answer, then $ f $ is called \emph{simultaneous Waring identifiable}. Identifiability is an extremely important property in the more general context of tensors, of which homogeneous polynomials are part, and it has applications in many different areas. For example, it is crucial in signal processing to solve the so called \emph{cocktail party problem}: in a room some people are talking, several receivers are set up to record all the speeches mixed together and the goal is to recover the voices of individuals. A similar situation is studied in chemistry, anytime it is necessary to \emph{unmix} the different chemicals in a solution. For other detailed applications, we refer to the book \cite{Lan}.  \\
\indent All the notions introduced so far are generalizations starting from the case of a single general polynomial, where Q1 and Q2  have been completely solved: the computation of the rank is due to Alexander-Hirschowitz and it can be found in \cite{AH}, while the list of all identifiable cases began with Sylvester (\cite{Sy}) and Hilbert (\cite{Hi}) and was recently completed in \cite{GM}. \\
\indent Concerning the case $ r>1 $, as far as I know, Q1 is still open while there are some partial results in the direction of Q2, in particular at the moment the known cases of simultaneous Waring identifiability are described in the following table:  

\begin{table}[h]
\begin{center}
\begin{tabular} {c |c |c |c |c }
{\bf n} & {\bf r} &  ${\bf d_{1}, \ldots, d_{r}}$ & {\bf k}  & {\bf Ref.} \\
$ 1 $ & $ \forall $ & $ d_1+1 \geq k  $  & $ \left \lceil \frac{1}{1+r} \sum_{j=1}^{r}{1+d_{j} \choose d_{j}} \right \rceil $ & \cite{CR} \\
$ \forall $ & $2$ & $2,2$ & $ n+1 $ & \cite{We} \\
$ 2 $ & $2$ & $2,3$ & $ 4 $ & \cite{Ro} \\
$ 2 $ & $3$ & $3,3,4$ & $ 7 $ & \cite{AGMO} \\
$ 2 $ & $4$ & $2,2,2,2$ & $ 4 $ & Veronese \\
\end{tabular}\vspace{0.2cm}
\end{center}
\end{table}

However it is not yet known whether the above list is complete or not. The aim of this paper is to deal with this classification problem, starting from the case $ n=2 $. More in details, let us focus on the three central rows of the previous table. In the recent paper \cite{BG} the authors show that, in the setting of pairs of ternary forms (i.e. $ n = r = 2 $) the unique identifiable examples are the case of two quadrics and the case of a quadric and a cubic. These cases were classically known, the first one goes back to Weierstrass (\cite{We}), while the second one to Roberts (\cite{Ro}). Notice that a general quadric has infinitely many decompositions with $ n+1 $ summands and a general ternary cubic of rank $ 4 $ is not identifiable, therefore by requiring simultaneous decompositions provides another useful technique to recover classical identifiability when it fails. In this direction, by adding the assumption that $ d_{2} - d_{1} = 1 $, the fact that Roberts' example is the unique identifiable one has already been proved in \cite{AGMO}. In the context of triples of ternary forms, the authors of \cite{AGMO} proved that a general polynomial vector with two cubics and a quartic, of rank $ 7 $, is identifiable. Roberts' example and the new case presented in \cite{AGMO} strongly suggested the following conjecture to obtain further simultaneous Waring identifiable instances: 
\begin{conj}\label{simconj}
For any $ d \geq 2 $, the general complex polynomial vector consisting of $ d-1 $ ternary forms of degree $ d $ and one ternary form of degree $ d+1 $ admits a unique simultaneous Waring decomposition with $ \frac{d^2+d+2}{2} $ summands.
\end{conj}
This conjecture arose during some fruitful discussions on the occasion of the annual cycle of seminars on \emph{Applied Algebraic Geometry}, that takes place approximately every three weeks from September to June between Firenze and Bologna and that involves many algebraic geometers, such as Alessandra Bernardi, Luca Chiantini, Alessandro Gimigliano, Massimiliano Mella, Giorgio Ottaviani together with Italian and foreign researchers, post-docs and phd students interested in this topic. In this sense, also the paper \cite{AGMO} represents a research product of one of the first series of seminars focusing on \emph{Numerical Algebraic Geometry}. \\
In this paper we show that Conjecture \ref{simconj} is not true, providing the counterexample described in the following:
\begin{thm}\label{counter}
The general complex polynomial vector $ f= (f_{1}, f_{2}, f_{3}, f_{4}) $, with $ f_{j} \in \C[x_{0}, x_{1}, x_{2}]_{4} $ for $ j \in \{1,2,3\} $ and $ f_{4} \in \C[x_{0}, x_{1}, x_{2}]_{5} $, has at least $ 36 $ simultaneous Waring decompositions with $ 11 $ summands.
\end{thm}   
The non-identifiability described in Theorem \ref{counter} has been obtained computationally by means of an algorithm implemented in \emph{Bertini} (\cite{Be}, \cite{BHSW}) and based on \emph{homotopy continuation} and \emph{monodromy loops} (\cite{HOOS}), in the spirit of \cite{AGMO} and \cite{A}. \\
\indent The paper is organized as follows. After introducing the subject, in Section $ 2 $ we fix main notations and definitions for the simultaneous Waring problem, together with the geometrical point of view. Finally, Section $ 3 $ is devoted to the computational procedure that allows us to prove Theorem \ref{counter} and then to detect the counter-example for Conjecture \ref{simconj}. \\
The ancillary file {\tt{decompositions.pdf}}, with the details on the $ 36 $ decompositions, completes the paper.

\section{Basics}

In this section we recall basic definitions and known results concerning the simultaneous Waring problem and its geometric translation via secant varieties. \\
\subsection{Simultaneous Waring background}\label{swb} $ \quad $\\

\indent In this paper we work over the complex field $ \C $. \\
Let $ (n, r) \in \mathbb{N}^2  - \{\underline{0}\} $ and let $ (d_{1}, \ldots, d_{r}) \in \mathbb{N}^r  - \{\underline{0}\} $ in ascending order. 
\begin{defn}\label{polvec}
A \emph{polynomial vector} $ f= (f_1, \ldots, f_r) $ is a vector of $ r $ complex forms in $ n+1 $ variables $ x_{0}, \ldots, x_{n} $ of degrees $ d_{1}, \ldots, d_{r} $, i.e. $ f_{j} \in \mathbb{C}[x_{0}, \ldots, x_{n}]_{d_{j}} $ for $ j \in \{1, \ldots, r\} $.  
\end{defn}
\begin{defn}\label{simWardec}
A complex \emph{simultaneous Waring decomposition} of $ f= (f_1, \ldots, f_r) $ is given by $ \ell_{1}, \ldots, \ell_{k} \in \mathbb{C}[x_{0}, \ldots, x_{n}]_{1} $ linear forms and $ (\lambda_{1}^{j}, \ldots, \lambda_{k}^{j}) \in \mathbb{C}^{k} - \{\underline{0}\} $ coefficients such that
\begin{equation}\label{vecWardec}
f = \sum_{i=1}^{k}\left(\lambda_{i}^{1}\ell_{i}^{d_{1}}, \ldots, \lambda_{i}^{r}\ell_{i}^{d_{r}}\right).
\end{equation}
\end{defn}
\begin{defn}\label{swr}
The complex \emph{simultaneous Waring rank} of a polynomial vector $ f = (f_{1}, \ldots, f_{r}) $  is the minimal positive integer $ k $ that appears in (\ref{vecWardec}).
\end{defn}
\noindent Notice that the notion introduced in Definition \ref{swr} yields a measure of the complexity of the polynomial vector $ f $. In particular, each $ \left(\lambda_{i}^{1}\ell_{i}^{d_{1}}, \ldots, \lambda_{i}^{r}\ell_{i}^{d_{r}}\right) $ has rank 1, being the simplest polynomial vector.
\begin{defn}
A complex rank-$ k $ polynomial vector $ f = (f_{1}, \ldots, f_{r}) $ is \emph{simultaneous Waring identifiable} if the expression (\ref{vecWardec}) is unique up to reordering the summands and rescaling the coefficients.
\end{defn}



\subsection{Secant varieties to projective bundles}\label{sec}$ \quad $\\

Let $\mathbb{G}(k,N)$ be the \emph{Grassmannian} of $ k-1 $-linear spaces in $ \mathbb{P}^{N-1} $, let $ X \subset \mathbb{P}^{N-1} $ be an irreducible variety and let 
$$ \Gamma_{k}(X) = \overline{\{((x_{1}, \ldots, x_{k}), L) \in \underbrace{X \times \ldots \times X}_{\text{k}} \times \mathbb{G}(k,N) \, | \, L = \langle x_{1}, \ldots, x_{k}\rangle\}}. $$
By construction, $ \Gamma_{k}(X) $ is irreducible and $ \dim\Gamma_{k}(X) = k \dim X $. Consider the projection onto the last factor $ \pi_{2}: \underbrace{X \times \ldots \times X}_{\text{k}} \times \mathbb{G}(k,N) \to \mathbb{G}(k,N) $ and denote $ \pi_{2}(\Gamma_{k}(X)) $ by $ S_{k}(X) $: $ S_{k}(X) $ is irreducible and $ \dim S_{k}(X) = \dim \Gamma_{k}(X) = k \dim X $. With the above notation, we give the following:
\begin{defn}
The \emph{abstract $k^{th}$-secant variety} of $ X $ is 
$$ sec_{k}(X) = \overline{\{(p,L) \in \mathbb{P}^{N-1} \times \mathbb{G}(k,N) \, | \, p \in L, L \in S_{k}(X)\}} $$
and the \emph{$k^{th}$-secant variety} of $ X $ is 
$$ \sigma_{k}(X) = \tilde{\pi}_{1}(sec_{k}(X)) $$
where $ \tilde{\pi}_{1} $ is the restriction to $ sec_{k}(X) $ of the projection $ \mathbb{P}^{N-1} \times \mathbb{G}(k,N) \to \mathbb{P}^{N-1} $.\\
$ X $ is \emph{$k$-defective} if 
$$ \dim \sigma_{k}(X) < \min\{\dim sec_{k}(X), N-1\} = \min\{k \dim X + k-1, N-1\}. $$
\end{defn}
\noindent A fundamental tool in the study of secant varieties is the following:
\begin{thm}[Terracini Lemma \cite{Te}\cite{CC}] \label{terracini}
Let $ X \subset \mathbb{P}^{N-1} $ be an irreducible variety and let $ z \in \sigma_{k}(X) $ be a general smooth point. The projective tangent space to $ \sigma_{k}(X) $ at $ z $ is given by
$$ \mathbb{T}_{z}\sigma_{k}(X) = \langle \mathbb{T}_{x_{1}}X,\ldots, \mathbb{T}_{x_{k}}X \rangle, $$
where $ x_{1}, \ldots, x_{k} $ are general smooth points in $ X $ such that $ z \in \langle x_{1}, \ldots, x_{k} \rangle $ and $ \mathbb{T}_{x_{i}}X $ denotes the projective tangent space to $ X $ at $ x_{i} $.  
\end{thm}
\begin{rem}
Notice that, if 
\begin{equation}\label{eq:equal}
\dim sec_{k}(X) = N-1,
\end{equation}
then the fiber of the map $ \tilde{\pi}_{1} $ is expected to be $ 0 $-dimensional, which implies that the general point $ p \in \mathbb{P}^{N-1} $ belongs to finitely many linear spaces cutting $ X $ in $ k $ points. Throughout this paper we work under the assumption that $ X $ is not $ k $-defective, so that the fiber of $ \tilde{\pi}_{1} $ is not empty. Therefore the condition (\ref{eq:equal}) implies that
\begin{equation}\label{eq:perfect}
k \dim X + k = N,
\end{equation}
which is called \emph{perfect case}. This is equivalent to say that $ \sigma_{k}(X) $ fills the ambient space $ \mathbb{P}^{N-1} $.
\end{rem}
In the simultaneous Waring setting, all the geometric tools introduced up to now are applied to the projective bundle $ X = \mathbb{P}(\mathcal{O}_{\mathbb{P}^{n}}(d_{1}) \oplus \ldots \oplus \mathcal{O}_{\mathbb{P}^{n}}(d_{r})) \subset \mathbb{P}^{N-1} =\mathbb{P}(H^{0}(\mathcal{O}_{\mathbb{P}^{n}}(d_{1})) \oplus \ldots, \oplus H^{0}(\mathcal{O}_{\mathbb{P}^{n}}(d_{r}))) $, with $ N = \sum_{i = 1}^{r} {n + d_{i} \choose d_{i}}$. For simplicity, throughout the paper we denote the projective bundle under investigation by $ X^{n}_{d_{1}, \ldots, d_{r}} $. For more details on $ X^{n}_{d_{1}, \ldots, d_{r}} $, we refer to \cite{AGMO}, where a detailed study on these varieties has been developed. In particular, regarding (\ref{vecWardec}), $ X^{n}_{d_{1}, \ldots, d_{r}} $ parametrizes the variety of rank-$ 1 $ polynomial vectors $ \left(\lambda_{i}^{1}\ell_{i}^{d_{1}}, \ldots, \lambda_{i}^{r}\ell_{i}^{d_{r}}\right) $, while $ f = (f_{1}, \ldots, f_{r}) $ belongs to $ \sigma_{k}(X^{n}_{d_{1}, \ldots, d_{r}}) $. Moreover, notice that $ \dim X^{n}_{d_{1}, \ldots, d_{r}} = n + r - 1 $, so that being in a perfect case, which represents the second assumption of our paper, means that 
\begin{equation}\label{eq:perfectsim}
k = {{\sum_{i = 1}^{r} {n + d_{i} \choose d_{i}}} \over {n+r}}.
\end{equation}
Therefore, under the assumptions that $ X^{n}_{d_{1}, \ldots, d_{r}} $ is not $ k $-defective and that its $ k $-secant variety fills the ambient space, the complex simultaneous Waring rank of a general polynomial vector $ f = (f_{1}, \ldots, f_{r}) $ is given in (\ref{eq:perfectsim}). This provides an answer to question Q1 introduced in \S\ref{sec:intr} under our two assumptions. More in general, since, in the simultaneous setting, the defectivity problem is still open, referring to Q1 there is the following conjecture, which goes in the direction of the main result in \cite{AH}: 
\begin{conj}
The complex simultaneous Waring rank of the general polynomial vector $ f = (f_{1}, \ldots, f_{r}) $, with parameters $ \{n, r, d_{1}, \ldots, d_{r}\} $ as in Definition \ref{polvec}, is 
$$ k = \left \lceil{ {{\sum_{i = 1}^{r} {n + d_{i} \choose d_{i}}} \over {n+r}} } \right \rceil,$$
where $ \left \lceil \, \right \rceil $ denotes the ceiling of the ratio.
\end{conj}
\begin{rem} 
For the purposes of this paper, that is to deal with Q2 and Conjecture \ref{simconj}, we focus on the projective bundle $ X^{2}_{\small{\underbrace{d, \ldots, d}_{\text{d-1}}},d+1} $, with $ d \geq 2 $. In this case, the positive integer $ k $ satisfying (\ref{eq:perfectsim}) is 
$$ {(d-1) {2 + d \choose d} + {3 + d \choose d+1} \over {2+d}} = {(d-1)(d+1)+(d+3) \over 2} = {d^2 + d + 2 \over 2}. $$
In particular, if $ d = 2 $, then one gets Roberts' example, while if $ d = 3 $, then we are dealing with the identifiable case discovered in \cite{AGMO}. Therefore, in order to check the truth of Conjecture \ref{simconj}, we start by considering the case $ d = 4 $ (and $ k = 11 $).  
\end{rem}

\section{Non-identifiability via Numerical Algebraic Geometry}

In this section we prove the main result of this paper. Indeed we describe the computational proof for Theorem \ref{counter}. \\

\noindent \emph{Proof of Theorem \ref{counter}}\\
\noindent The first step is to check that the projective bundle $ X^{2}_{4,4,4,5} $ is not $11$-defective and that $ \dim \sigma_{11}(X^{2}_{4,4,4,5}) = 65 $, which can be done, in the spirit of Theorem \ref{terracini}, by computing the dimension of the span of the affine tangent spaces to $ X^{2}_{4,4,4,5} $ at $ 11 $ (rational) random points and then by applying semi-continuity properties. In the following, the script of the probabilistic algorithm, implemented in Macaulay2 software \cite{M2}, that allows us to check the truth of the required hypothesis:

{\small{\begin{verbatim}
R = QQ[a_0,a_1,a_2,x,y,z];
aa = matrix{{a_0,a_1,a_2}};
a4 = symmetricPower(4,aa);
a5 = symmetricPower(5,aa);
j = diff(transpose basis(1,R),a4|x*a4|y*a4|z*a5);
j1 = sub(j,random(R^{1:0},R^{6:0})); -- matrix 6*66
for i from 2 to 11 do j1 = j1||sub(j,random(R^{1:0},R^{6:0})) -- matrix 66*66
rank j1 -- 66
\end{verbatim}}}

Then we solve, via the software Bertini for \emph{Numerical Algebraic Geometry} \cite{Be}, the polynomial system (\ref{vecWardec}), where $ n=2, r=4, d_{1}=d_{2}=d_{3}=4, d_{4} = 5, k=11 $:

\begin{equation}\label{eq:polsys}
\left\{
\begin{array}{l}
f_{1} - \lambda_{1}^{1}\ell_{1}^{4}- \ldots -  \lambda_{11}^{1}\ell_{11}^{4} = 0 \\
f_{2} - \lambda_{1}^{2}\ell_{1}^{4}- \ldots -  \lambda_{11}^{2}\ell_{11}^{4} = 0 \\
f_{3} - \lambda_{1}^{3}\ell_{1}^{4}- \ldots -  \lambda_{11}^{3}\ell_{11}^{4} = 0 \\
f_{4} - \lambda_{1}^{4}\ell_{1}^{5}- \ldots -  \lambda_{11}^{4}\ell_{11}^{5} = 0 \\
\end{array}
\right. .
\end{equation}

\noindent Notice that $ f_{1}, f_{2}, f_{3} \in \mathbb{C}[x_{0}, x_{1}, x_{2}]_{4} $ and $  f_{4} \in \mathbb{C}[x_{0}, x_{1}, x_{2}]_{5} $ are known, while $ \ell_{i} = x_{0}+ l_{1}^{i}x_{1}+ l_{2}^{i}x_{2} \in \C[x_{0}, x_{1}, x_{2}]_{1}$ and $ \lambda_{i}^{j} \in \C $ are unknown, for $ i \in \{1, \ldots, 11\} $ and $ j \in \{1,2,3,4\} $. By developing computations, each of the first three equations of (\ref{eq:polsys}) provides $ 15 $ conditions, while the fourth equation gives $ 21 $ conditions. Therefore (\ref{eq:polsys}) is equivalent to a non linear system consisting of $ 66 $ equations and unknowns, usually denoted by $ F_{\left(f_{1}, f_{2}, f_{3}, f_{4}\right)}\left(\left[l_{1}^{1}, l_{2}^{1}, \lambda_{1}^{1}, \lambda_{1}^{2},\lambda_{1}^{3},\lambda_{1}^{4}\right],\ldots, \left[l_{1}^{11}, l_{2}^{11}, \lambda_{11}^{1}, \lambda_{11}^{2},\lambda_{11}^{3},\lambda_{11}^{4}\right]\right) $, to give evidence of parameters and unknowns. The equations of (\ref{eq:polsys}) can be constructed by means of Macaulay2, for example via the following script:

{\small{\begin{verbatim}
k = ZZ/101;
S = k[v_0..v_65,x_0,x_1,x_2];
xx = matrix{{x_0,x_1,x_2}};
f1 = diff(symmetricPower(4,xx),sum(11,i->v_(6*i+2)*(matrix{{1,v_(6*i),
     v_(6*i+1)}}*transpose xx)^4)) -- 15 components
f2 = diff(symmetricPower(4,xx),sum(11,i->v_(6*i+3)*(matrix{{1,v_(6*i),
     v_(6*i+1)}}*transpose xx)^4)) -- 15 components
f3 = diff(symmetricPower(4,xx),sum(11,i->v_(6*i+4)*(matrix{{1,v_(6*i),
     v_(6*i+1)}}*transpose xx)^4)) -- 15 components
f4 = diff(symmetricPower(5,xx),sum(11,i->v_(6*i+5)*(matrix{{1,v_(6*i),
     v_(6*i+1)}}*transpose xx)^5)) -- 21 components
\end{verbatim}}}
 
\noindent The solutions of $ F_{\left(f_{1}, f_{2}, f_{3}, f_{4}\right)}\left(\left[l_{1}^{1}, l_{2}^{1}, \lambda_{1}^{1}, \lambda_{1}^{2},\lambda_{1}^{3},\lambda_{1}^{4}\right],\ldots, \left[l_{1}^{11}, l_{2}^{11}, \lambda_{11}^{1}, \lambda_{11}^{2},\lambda_{11}^{3},\lambda_{11}^{4}\right]\right) $, seen as vectors with $ 66 $ components, provide the complex simultaneous Waring decompositions of length $ 11 $ of the polynomial vector $ f = (f_{1}, f_{2}, f_{3}, f_{4}) $. We aim to show that, for a general $ f $, the number of different solutions is $ 36 $ instead of $ 1 $. From this point on, all the scripts of the procedure are contained in the ancillary folder {\tt{ternarie\_4445}}. \\
\noindent Since $ f $ needs to be general, we first assign random complex entries to the vector of unknowns, obtaining {\small{$ \left(\left[\overline{l}_{1}^{1}, \overline{l}_{2}^{1}, \overline{\lambda}_{1}^{1},\overline{\lambda}_{1}^{2},\overline{\lambda}_{1}^{3}, \overline{\lambda}_{1}^{4}\right], \ldots, \left[\overline{l}_{1}^{11}, \overline{l}_{2}^{11}, \overline{\lambda}_{11}^{1}, \overline{\lambda}_{11}^{2},\overline{\lambda}_{11}^{3},\overline{\lambda}_{11}^{4}\right]\right) $}}, the \emph{start-point} of the procedure. By substituting these values to the unknowns of $ F_{\left(f_{1}, f_{2}, f_{3}, f_{4}\right)} $, we get a general polynomial vector $ \overline{f} = (\overline{f}_{1}, \overline{f}_{2},\overline{f}_{3}, \overline{f}_{4}) $, which $ 66 $ coefficients are the \emph{start-parameters} of our procedure. This is done by means of the scripts {\tt{input\_eval}} and {\tt{start\_good}} appearing in the sub-folder {\tt{Eval}}, where the latter file contains the start-point and the former the 66 equations that need to be evaluated at the coordinates of the start-point (in particular, in {\tt{input\_eval}}, $ v_{0}, \ldots, v_{65} $ are the variables to which the random entries need to be assigned while $ f_{0}, \ldots, f_{65} $ denote the equations of the polynomial system; the {\tt{TrackType}} configuration inside {\tt{input\_eval}} is $ -4 $ since at this step Bertini evaluates the input system at the given points). The start-parameters are obtained with this evaluation and are contained in the file {\tt{start\_parameters\_good}} of the sub-folder {\tt{Monodromy}}. Our goal is to prove that $ \overline{f} $ has $ 35 $ other decompositions, in addition to the one given by the start-point. \\
We constructed our counter-example by means of the following start-point, which we randomly generated:

{\small{$$\left[\overline{l}_{1}^{1}, \overline{l}_{2}^{1}, \overline{\lambda}_{1}^{1},\overline{\lambda}_{1}^{2},\overline{\lambda}_{1}^{3}, \overline{\lambda}_{1}^{4}\right] = \left[-9.178673081212201\cdot10^{-2}- i \, 5.245599426621500\cdot10^{-2}, \right. $$
$$ \quad\quad\quad\quad\quad\quad\quad\quad\quad\quad\, -8.825275325577600\cdot10^{-2}- i \, 2.155292826256500\cdot10^{-2},  $$
$$ \quad\quad\quad\quad\quad\quad\quad\quad\quad\quad\quad 8.022566436298600\cdot10^{-1} + i \, 9.927688077008801\cdot10^{-2}, $$
$$ \quad\quad\quad\quad\quad\quad\quad\quad\quad\quad\quad 7.524729028260101\cdot10^{-2} + i \, 8.367569801775700\cdot10^{-2}, $$
$$  \quad\quad\quad\quad\quad\quad\quad\quad\quad\quad -1.250209354876000\cdot10^{-3} - i \, 6.530239085633700\cdot10^{-2}, $$
$$ \quad\quad\quad\quad\quad\quad\quad\quad\quad\quad\quad \left. 4.999653700389400\cdot10^{-2} + i \, 9.303029998119500\cdot10^{-2} \right] $$

$$\left[\overline{l}_{1}^{2}, \overline{l}_{2}^{2}, \overline{\lambda}_{2}^{1},\overline{\lambda}_{2}^{2},\overline{\lambda}_{2}^{3}, \overline{\lambda}_{2}^{4}\right] = \left[1.195416802244310\cdot10^{-2}- i \, 2.225934180354000\cdot10^{-3}, \right. $$
$$ \quad\quad\quad\quad\quad\quad\quad\quad\quad\quad\,\,\, -6.560577710833700\cdot10^{-1} -  i \, 9.891546819700110\cdot10^{-2},  $$
$$ \quad\quad\quad\quad\quad\quad\quad\quad\quad\quad\,\,\, -4.115491104918700\cdot10^{-2} + i \, 4.092292289326200\cdot10^{-2}, $$
$$ \quad\quad\quad\quad\quad\quad\quad\quad\quad\quad\,\,\, -2.570701095199700\cdot10^{-2} - i \, 8.005816221966300\cdot10^{-2}, $$
$$  \quad\quad\quad\quad\quad\quad\quad\quad\quad\quad 2.941914008262100\cdot10^{-2} - i \, 7.332772392640199\cdot10^{-2}, $$
$$ \quad\quad\quad\quad\quad\quad\quad\quad\quad\quad\,\,\, \left. -7.507007824205600\cdot10^{-2} + i \, 2.699645975995900\cdot10^{-2} \right] $$

$$\left[\overline{l}_{1}^{3}, \overline{l}_{2}^{3}, \overline{\lambda}_{3}^{1},\overline{\lambda}_{3}^{2},\overline{\lambda}_{3}^{3}, \overline{\lambda}_{3}^{4}\right] = \left[8.306470571337600\cdot10^{-2}+ i \, 6.925123923475100\cdot10^{-2}, \right. $$
$$ \quad\quad\quad\quad\quad\quad\quad\quad\quad\quad\,\,\, -2.631866554121000\cdot10^{-3} +  i \, 6.450055333951101\cdot10^{-2},  $$
$$ \quad\quad\quad\quad\quad\quad\quad\quad\quad\quad\,\,\, -2.520169416138000\cdot10^{-2} + i \, 7.947572402054000\cdot10^{-2}, $$
$$ \quad\quad\quad\quad\quad\quad\quad\quad\quad\quad 6.628901260347200\cdot10^{-2} - i \, 9.954313791954100\cdot10^{-2}, $$
$$  \quad\quad\quad\quad\quad\quad\quad\quad\quad\quad\,\, -4.986358377275200\cdot10^{-2} - i \, 6.159247708199100\cdot10^{-2}, $$
$$ \quad\quad\quad\quad\quad\quad\quad\quad\quad\quad\,\,\, \left. -4.363786207688800\cdot10^{-2} + i \, 6.858262771755810\cdot10^{-2} \right] $$

$$\left[\overline{l}_{1}^{4}, \overline{l}_{2}^{4}, \overline{\lambda}_{4}^{1},\overline{\lambda}_{4}^{2},\overline{\lambda}_{4}^{3}, \overline{\lambda}_{4}^{4}\right] = \left[-3.777699832746000\cdot10^{-2} - i \, 9.873023399290000\cdot10^{-2}, \right. $$
$$ \quad\quad\quad\quad\quad\quad\quad\quad\quad\,\,\,\, 2.251959174854600\cdot10^{-2} -  i \, 9.135044625682701\cdot10^{-2},  $$
$$ \quad\quad\quad\quad\quad\quad\quad\quad\quad\,\,\,\,\, 3.669676250951100\cdot10^{-2} + i \, 7.622324724014400\cdot10^{-2}, $$
$$ \quad\quad\quad\quad\quad\quad\quad\quad\quad\,\,\,\,\, 3.634917018379500\cdot10^{-2} - i \, 5.107012333167600\cdot10^{-2}, $$
$$  \quad\quad\quad\quad\quad\quad\quad\quad\quad\quad\,\, -7.505713542195000\cdot10^{-3} + i \,  5.154623436746600\cdot10^{-2}, $$
$$ \quad\quad\quad\quad\quad\quad\quad\quad\quad\quad\,\,\, \left. -5.543948948927200\cdot10^{-2} + i \, 9.026674172056499\cdot10^{-2} \right] $$

$$\left[\overline{l}_{1}^{5}, \overline{l}_{2}^{5}, \overline{\lambda}_{5}^{1},\overline{\lambda}_{5}^{2},\overline{\lambda}_{5}^{3}, \overline{\lambda}_{5}^{4}\right] = \left[9.015391063311000\cdot10^{-2} + i \, 2.074965436740300\cdot10^{-2}, \right. $$
$$ \quad\quad\quad\quad\quad\quad\quad\quad\quad\quad\, 2.113776820185300\cdot10^{-2} +  i \, 4.703248826760600\cdot10^{-2},  $$
$$ \quad\quad\quad\quad\quad\quad\quad\quad\quad\quad\, 6.374079255711100\cdot10^{-2} + i \, 8.650339728394700\cdot10^{-2}, $$
$$ \quad\quad\quad\quad\quad\quad\quad\quad\quad\quad\, 4.048211909924600\cdot10^{-2} - i \, 7.207606064662800\cdot10^{-2}, $$
$$  \quad\quad\quad\quad\quad\quad\quad\quad\quad\quad\,\,\, -8.151521021270800\cdot10^{-2} - i \, 2.388618850411100\cdot10^{-2}, $$
$$ \quad\quad\quad\quad\quad\quad\quad\quad\quad\quad\,\,\, \left. -4.363786207688800\cdot10^{-2} + i \, 6.858262771755810\cdot10^{-2} \right] $$

$$\left[\overline{l}_{1}^{6}, \overline{l}_{2}^{6}, \overline{\lambda}_{6}^{1},\overline{\lambda}_{6}^{2},\overline{\lambda}_{6}^{3}, \overline{\lambda}_{6}^{4}\right] = \left[-3.777699832746000\cdot10^{-2} - i \, 9.873023399290000\cdot10^{-2}, \right. $$
$$ \quad\quad\quad\quad\quad\quad\quad\quad\quad\quad\, 2.251959174854600\cdot10^{-2} -  i \, 1.135044625682700\cdot10^{-2},  $$
$$ \quad\quad\quad\quad\quad\quad\quad\quad\quad\quad\, 3.669676250951100\cdot10^{-2} + i \, 3.622324724014400\cdot10^{-2}, $$
$$ \quad\quad\quad\quad\quad\quad\quad\quad\quad\quad\, 3.634917018379500\cdot10^{-2} - i \, 5.107012333167600\cdot10^{-2}, $$
$$  \quad\quad\quad\quad\quad\quad\quad\quad\quad\quad\,\,\, -7.505713542195000\cdot10^{-3} + i \, 1.154623436746600\cdot10^{-2}, $$
$$ \quad\quad\quad\quad\quad\quad\quad\quad\quad\quad\,\,\, \left. -5.543948948927200\cdot10^{-2} + i \, 9.026674172056499\cdot10^{-2} \right] $$

$$\left[\overline{l}_{1}^{7}, \overline{l}_{2}^{7}, \overline{\lambda}_{7}^{1},\overline{\lambda}_{7}^{2},\overline{\lambda}_{7}^{3}, \overline{\lambda}_{7}^{4}\right] = \left[9.015391063311000\cdot10^{-2} + i \, 2.074965436740300\cdot10^{-2}, \right. $$
$$ \quad\quad\quad\quad\quad\quad\quad\quad\quad\quad\, 1.137768201853000\cdot10^{-3} +  i \, 4.703248826760600\cdot10^{-2},  $$
$$ \quad\quad\quad\quad\quad\quad\quad\quad\quad\quad\, 6.374079255711100\cdot10^{-2} + i \, 8.650339728394700\cdot10^{-2}, $$
$$ \quad\quad\quad\quad\quad\quad\quad\quad\quad\quad\, 4.948211909924600\cdot10^{-2} - i \, 7.207606064662800\cdot10^{-2}, $$
$$  \quad\quad\quad\quad\quad\quad\quad\quad\quad\quad\,\,\,\,\, -8.151521021270800\cdot10^{-2} - i \, 2.388618850411100\cdot10^{-2}, $$
$$ \quad\quad\quad\quad\quad\quad\quad\quad\quad\quad\, \left. 4.177457890165700\cdot10^{-2} - i \, 6.689333670048300\cdot10^{-2} \right] $$

$$\left[\overline{l}_{1}^{8}, \overline{l}_{2}^{8}, \overline{\lambda}_{8}^{1},\overline{\lambda}_{8}^{2},\overline{\lambda}_{8}^{3}, \overline{\lambda}_{8}^{4}\right] = \left[-2.952938879084900\cdot10^{-2} - i \, 4.829159279506610\cdot10^{-2}, \right. $$
$$ \quad\quad\quad\quad\quad\quad\quad\quad\quad\quad\, -6.105800366706000\cdot10^{-3} -  i \, 9.111150804119310\cdot10^{-2},  $$
$$ \quad\quad\quad\quad\quad\quad\quad\quad\quad\quad\, -4.424765849460100\cdot10^{-3} - i \, 5.448494845808000\cdot10^{-3}, $$
$$ \quad\quad\quad\quad\quad\quad\quad\quad\quad\quad\, -8.459618609850000\cdot10^{-4} + i \, 9.105717695359500\cdot10^{-2}, $$
$$  \quad\quad\quad\quad\quad\quad\quad\quad\quad\quad 9.212993622611000\cdot10^{-3} + i \, 7.622100043006000\cdot10^{-3}, $$
$$ \quad\quad\quad\quad\quad\quad\quad\quad\quad\quad \left. 4.814038260778000\cdot10^{-2} + i \, 3.339139584718200\cdot10^{-2} \right] $$

$$\left[\overline{l}_{1}^{9}, \overline{l}_{2}^{9}, \overline{\lambda}_{9}^{1},\overline{\lambda}_{9}^{2},\overline{\lambda}_{9}^{3}, \overline{\lambda}_{9}^{4}\right] = \left[3.333163313138000\cdot10^{-2} + i \, 8.987872802854000\cdot10^{-2}, \right. $$
$$ \quad\quad\quad\quad\quad\quad\quad\quad\quad\quad\, 4.924941240444200\cdot10^{-2} +  i \, 4.884444793924100\cdot10^{-2},  $$
$$ \quad\quad\quad\quad\quad\quad\quad\quad\quad\quad\,\,\, -4.344354377175200\cdot10^{-2} - i \, 2.259235703299100\cdot10^{-2}, $$
$$ \quad\quad\quad\quad\quad\quad\quad\quad\quad\quad\,\,\, -7.463782067688800\cdot10^{-2} + i \, 7.758622222255809\cdot10^{-2}, $$
$$  \quad\quad\quad\quad\quad\quad\quad\quad\quad\quad\,\,\, -7.994595436034500\cdot10^{-2} + i \, 5.510716975355500\cdot10^{-2}, $$
$$ \quad\quad\quad\quad\quad\quad\quad\quad\quad\quad \left. 8.955399362261099\cdot10^{-1} + i \, 9.456553464300600\cdot10^{-2} \right] $$

$$\left[\overline{l}_{1}^{10}, \overline{l}_{2}^{10}, \overline{\lambda}_{10}^{1},\overline{\lambda}_{10}^{2},\overline{\lambda}_{10}^{3}, \overline{\lambda}_{10}^{4}\right] = \left[-5.32316331313800\cdot10^{-2} + i \, 3.337872802854000\cdot10^{-2}, \right. $$
$$ \quad\quad\quad\quad\quad\quad\quad\quad\quad\quad\,\,\,\,\,\,\,\,\,\,\,\,\,\,\,\, 3.334941240444200\cdot10^{-2} +  i \, 4.444444793924100\cdot10^{-2},  $$
$$ \quad\quad\quad\quad\quad\quad\quad\quad\quad\quad\,\,\,\,\,\,\,\,\,\,\,\,\,\,\,\,\,\,\,\,\, -4.333354377175200\cdot10^{-2} - i \, -2.259235703299100\cdot10^{-2}, $$
$$ \quad\quad\quad\quad\quad\quad\quad\quad\quad\quad\,\,\,\,\,\,\,\,\,\,\,\,\,\,\,\,\,\, -9.763782067688800\cdot10^{-2} + i \, 6.586222222558100\cdot10^{-3}, $$
$$  \quad\quad\quad\quad\quad\quad\quad\quad\quad\quad\,\,\,\,\,\,\,\,\,\,\,\,\,\,\,\,\,\, -4.445954360345000\cdot10^{-3} + i \, 9.510716975355500\cdot10^{-2}, $$
$$ \quad\quad\quad\quad\quad\quad\quad\quad\quad\quad\,\,\,\,\,\,\,\,\,\,\,\,\,\, \left. 8.034399362261100\cdot10^{-1} + i \, 9.056213464300600\cdot10^{-2} \right] $$

$$\left[\overline{l}_{1}^{11}, \overline{l}_{2}^{11}, \overline{\lambda}_{11}^{1},\overline{\lambda}_{11}^{2},\overline{\lambda}_{11}^{3}, \overline{\lambda}_{11}^{4}\right] = \left[5.978211909924600\cdot10^{-2} + i \, 1.237606064662800\cdot10^{-2}, \right. $$
$$ \quad\quad\quad\quad\quad\quad\quad\quad\quad\quad\,\,\,\,\,\,\,\,\,\,\,\,\,\,\,\, 8.961521021270801\cdot10^{-2} +  i \, 3.277618850411100\cdot10^{-2},  $$
$$ \quad\quad\quad\quad\quad\quad\quad\quad\quad\quad\,\,\,\,\,\,\,\,\,\,\,\,\,\,\,\,\,\,\,\,\, -4.377457890165700\cdot10^{-2} + i \, 5.579333670048300\cdot10^{-2}, $$
$$ \quad\quad\quad\quad\quad\quad\quad\quad\quad\quad\,\,\,\,\,\,\,\,\,\,\,\,\,\,\, 9.252938879084900\cdot10^{-2} - i \, 4.845159279506610\cdot10^{-2}, $$
$$  \quad\quad\quad\quad\quad\quad\quad\quad\quad\quad\,\,\,\,\,\,\,\,\,\,\,\,\,\, 3.610580036670600\cdot10^{-2} + i \, 1.922150804119310\cdot10^{-2}, $$
$$ \quad\quad\quad\quad\quad\quad\quad\quad\quad\quad\,\,\,\,\,\,\,\,\,\,\,\,\,\,\,\,\, \left. 5.177457890165700\cdot10^{-2} + i \,  5.571333670048300\cdot10^{-2} \right]. $$}}

\noindent Therefore we consider the polynomial system $ F_{\left(\overline{f}_{1},\overline{f}_{2},\overline{f}_{3}, \overline{f}_{4}\right)} $, of which we know a solution, {\small{$ \left(\left[\overline{l}_{1}^{1}, \overline{l}_{2}^{1}, \overline{\lambda}_{1}^{1},\overline{\lambda}_{1}^{2},\overline{\lambda}_{1}^{3}, \overline{\lambda}_{1}^{4}\right], \ldots, \left[\overline{l}_{1}^{11}, \overline{l}_{2}^{11}, \overline{\lambda}_{11}^{1}, \overline{\lambda}_{11}^{2},\overline{\lambda}_{11}^{3},\overline{\lambda}_{11}^{4}\right]\right) $}}, and we start a monodromy loop divided in $ 3 $ linear steps. \\
\noindent In the first step, we change the start-parameters of $ F_{\left(\overline{f}_{1},\overline{f}_{2},\overline{f}_{3}, \overline{f}_{4}\right)} $ with random complex numbers (\emph{final-parameters}) getting a different polynomial system $ F_{1} $, still consisting of $ 66 $ equations and unknowns. By means of a \emph{segment homotopy} $ H_{0} : \C^{66} \times [0,1] \to \C^{66} $
between $ F_{\left(\overline{f}_{1}, \overline{f}_{2},\overline{f}_{3}, \overline{f}_{4}\right)} $ and $ F_{1} $ we thus connect the start-point {\small{$ \left(\left[\overline{l}_{1}^{1}, \overline{l}_{2}^{1}, \overline{\lambda}_{1}^{1},\overline{\lambda}_{1}^{2},\overline{\lambda}_{1}^{3}, \overline{\lambda}_{1}^{4}\right], \ldots, \left[\overline{l}_{1}^{11}, \overline{l}_{2}^{11}, \overline{\lambda}_{11}^{1}, \overline{\lambda}_{11}^{2},\overline{\lambda}_{11}^{3},\overline{\lambda}_{11}^{4}\right]\right) $}} to a solution (\emph{end-point}) of $ F_{1} $. \\
\noindent Then we replace $ F_{\left(\overline{f}_{1},\overline{f}_{2},\overline{f}_{3}, \overline{f}_{4}\right)} $, {\small{$ \left(\left[\overline{l}_{1}^{1}, \overline{l}_{2}^{1}, \overline{\lambda}_{1}^{1},\overline{\lambda}_{1}^{2},\overline{\lambda}_{1}^{3}, \overline{\lambda}_{1}^{4}\right], \ldots, \left[\overline{l}_{1}^{11}, \overline{l}_{2}^{11}, \overline{\lambda}_{11}^{1}, \overline{\lambda}_{11}^{2},\overline{\lambda}_{11}^{3},\overline{\lambda}_{11}^{4}\right]\right) $}} and $ \left(\overline{f}_{1},\overline{f}_{2},\overline{f}_{3}, \overline{f}_{4}\right) $ with $ F_{1} $, the end-point and the final parameters and we iterate the procedure of the first step. We thus obtain a polynomial system $ F_{2} $ and a segment homotopy $ H_{1} : \C^{66} \times [0,1] \to \C^{66} $ linking a solution of $ F_{1} $ (the output of the first step) to a solution of $ F_{2} $. \\
\noindent We finally close the loop by means of a segment homotopy $ H_{2} : \C^{66} \times [0,1] \to \C^{66} $ which sends a solution of $ F_{2} $ (the output of the second step) to a solution of $ F_{\left(\overline{f}_{1}, \overline{f}_{2},\overline{f}_{3}, \overline{f}_{4}\right)} $. \\
We point out that the equations of $ F_{1} $, $ F_{2} $ and $ F_{\left(\overline{f}_{1}, \overline{f}_{2},\overline{f}_{3}, \overline{f}_{4}\right)} $ (seen as the three \emph{vertexes} of a triangle loop) are obtained by changing in the file {\tt{input\_move}} (contained in the sub-folder {\tt{Monodromy}}) the parameters $p_{0}, \ldots, p_{65}$ with random entries; the segment homotopies $ H_{0}, H_{1}, H_{2} $ are constructed by Bertini (this explains why in {\tt{input\_move}} the instruction {\tt{ParameterHomotopy}} is set on the value $2$). By means of the file {\tt{input\_Newton}} (of the sub-folder {\tt{Monodromy}}), Bertini performs Newton's method (which explains why {\tt{TrackType}} is set on the value $ -2 $) in order to sharpen the solutions obtained at the end of a triangle loop.  \\
\noindent Since the output of the loop is not {\small{$ \left(\left[\overline{l}_{1}^{1}, \overline{l}_{2}^{1}, \overline{\lambda}_{1}^{1},\overline{\lambda}_{1}^{2},\overline{\lambda}_{1}^{3}, \overline{\lambda}_{1}^{4}\right], \ldots, \left[\overline{l}_{1}^{6}, \overline{l}_{2}^{6}, \overline{\lambda}_{6}^{1}, \overline{\lambda}_{6}^{2},\overline{\lambda}_{6}^{3},\overline{\lambda}_{6}^{4}\right]\right)$}}, then $ \overline{f} = (\overline{f}_{1}, \overline{f}_{2},\overline{f}_{3}, \overline{f}_{4}) $ provides a counter-example to Conjecture \ref{simconj}. In particular, by iterating the procedure many times, we get that the number of different solutions of $ F_{\left(\overline{f}_{1}, \overline{f}_{2},\overline{f}_{3}, \overline{f}_{4}\right)} $ stabilizes on $ 36 $, which concludes the proof of Theorem \ref{counter}. \\

\noindent More details on the $ 36 $ decompositions can be found in the ancillary file\\
{\tt{decompositions.pdf}}.  \\

\noindent\emph{Acknowledgments}. I would like to thank the anonymous referees for their valuable suggestions and useful comments which contributed to improve this paper. \\

\end{document}